\documentclass[11pt]{amsart}

\usepackage{latexsym,amsthm,amssymb}

\usepackage[leqno]{amsmath}

\usepackage{mathrsfs}

\usepackage[all]{xy}

\newtheorem{theorem}{Theorem}[section]
\newtheorem{lemma}[theorem]{Lemma}
\newtheorem{definition}[theorem]{Definition}
\newtheorem{defi}[theorem]{Definition}

\newtheorem{proposition}[theorem]{Proposition}
\newtheorem{cor}[theorem]{Corollary}
\newtheorem{remark}[theorem]{Remark}

\numberwithin{equation}{theorem}

\newenvironment{Proof}
{\noindent \emph{Proof.}}
{\hfill $\square$ \medskip}

\renewcommand{\mathcal}{\mathscr}
\newcommand{\Cal}{\mathcal}

\newcommand{\SC}{{\mathcal{C}}}
\newcommand{\SD}{{\mathcal{D}}}
\newcommand{\SE}{{\mathcal{E}}}
\newcommand{\SF}{{\mathcal{F}}}

\newcommand{\SI}{{\mathcal{I}}}

\newcommand{\SL}{{\mathcal{L}}}
\newcommand{\SM}{{\mathcal{M}}}
\newcommand{\SN}{{\mathcal{N}}}
\newcommand{\SO}{{\mathcal{O}}}
\newcommand{\SP}{{\mathcal{P}}}

\newcommand{\SX}{{\mathcal{X}}}
\newcommand{\SY}{{\mathcal{Y}}}

\newcommand{\PP}{\bold{P}}

\newcommand{\WY}{\widetilde{Y}}
\newcommand{\WX}{\widetilde{X}}

\newcommand{\WL}{\widetilde{L}}

\newcommand{\wi}{\widetilde{i}}

\newcommand{\wphi}{\widetilde{\varphi}}
\newcommand{\wrho}{\widetilde{\rho}}
\newcommand{\WPsi}{\widetilde{\Psi}}

\newcommand{\womega}{\widetilde{\omega}}

\title[Deformation of finite morphisms and smoothing
of ropes]{Deformation of finite morphisms and
\\smoothing of ropes}

\author{Francisco Javier Gallego}
\author{Miguel Gonz\'alez}
\author{\\ Bangere P. Purnaprajna}

\address{Departamento de \'Algebra, Universidad Complutense de Madrid}
\email{gallego@mat.ucm.es}
\address{Departamento de \'Algebra, Universidad Complutense de Madrid}
\email{mgonza@mat.ucm.es}
\address{Department of Mathematics, University of Kansas}
\email{purna@math.ku.edu}

\begin{document}
\begin{abstract}
In this article we prove that most ropes of arbitrary multiplicity,
supported on smooth curves can be smoothed. By
a rope being smoothable we
mean that the rope is the flat limit of a family of smooth,
irreducible curves. To construct a smoothing,
we connect, on the one hand,
deformations of a finite morphism to projective space
and, on the other hand, morphisms from a rope to projective space.
We  also prove a general result of independent interest, namely, finite
covers onto smooth
irreducible curves embedded in projective space can be deformed
to a family of $1:1$ maps.
We apply our general
theory to prove the smoothing of ropes of multiplicity $3$ on
$\bold P^1$.
Even though this article focuses on ropes of dimension $1$, our method
yields a general approach to deal with the smoothing of ropes of
higher dimension.
\end{abstract}
\maketitle

\section{Introduction}\label{intro}

This article contains two themes. The first part presents a general
method to deal with ropes, their smoothings and their relation to
deformations of finite morphisms. This method
produces, among other things, two interesting, not previously known
results.
The first one is the fact that finite
covers onto smooth
irreducible curves embedded in projective space can be deformed
to a family of $1:1$ maps. The second one is the smoothing of most ropes of
arbitrary multiplicity,
supported on smooth curves.
In the second part of the article we carry out a
detailed study of the projective embeddings of ropes of multiplicity
$3$ on $\bold P^1$ and we apply the general theory developed in the
first part to prove that they are smoothable.

\medskip

Let $Y$ be a smooth, irreducible projective curve. A rope $\WY$ on $Y$
of multiplicity $m$ is a nowhere reduced scheme whose reduced
structure is
$Y$ and which locally looks like  the first infinitesimal neighborhood
of
$Y$ inside the total space of a vector bundle of rank $m-1$. Such a
scheme
is projective and, even though singular everywhere, it is locally
Cohen-Macaulay. The ideal sheaf of $Y$ inside of $\WY$ happens to be a
locally free sheaf $\Cal E$ of rank $m-1$, the so-called conormal
bundle
of $\WY$. The conormal bundle $\Cal E$ is an important invariant of
the
rope and, together with $Y$, determines the arithmetic genus of
$\WY$. On
the other hand, fixed $Y$ and $\Cal E$, there are in general many
non--isomorphic rope structures $\WY$ on $Y$ having conormal bundle
$\Cal E$. Indeed, these ropes structures are parametrized by
Ext$^1(\Omega_Y, \Cal E)=H^1(\Omega^*_Y \otimes \Cal E)$. For an
introduction on ropes, the reader can look into~\cite{BE},
where Bayer and Eisenbud study at length ropes of multiplicity $2$,
the so-called ribbons.

\medskip

Multiple structures appear in many contexts in Algebraic Geometry.
For instance, the ubiquitously used so-called Ferrand--Szpiro doubling
(that is, multiplicity two structures) is a very special case of the
structures of arbitrary multiplicity dealt with in this article.
The Ferrand--Szpiro doubling is
used, mostly in an affine setting,
in the theory of projective modules and in problems dealing
with complete intersections, among other places. Projective multiple
structures, which are more complicated, are the ones studied in this paper.
The so-called canonical ribbons are another special case of these
multiple
structures
and attracted considerable attention during the 90s.
\medskip

From a classical point of view, a geometer is primarily interested in
smooth objects, or, at any rate, varieties. Thus, there are two
questions one might
reasonably ask, namely, whether ropes arise in a natural, uncontrived
geometrical context and what the connection of ropes to smooth
varieties is.
The answer to both matters is satisfactory. Precisely,
(see Theorems~\ref{coversmoothing} and~\ref{embsmoothing} and Corollary ~\ref{absgensmooth}), a rope $\WY$
with nonnegative arithmetic genus always arises as a flat limit of the
images of a family of projective
embeddings of smooth, irreducible curves, when those embeddings degenerate to
$\pi$, provided that $\Cal E$ satisfies the following condition:
%
%\smallskip
%
%
%
$$\displaylines{(*) \quad \Cal E \text{ can be realized as the trace
    zero module of} \cr
\text{a
smooth, irreducible,
finite cover } X \overset{\pi} \longrightarrow Y.}$$

%\smallskip

\noindent This answers right away the first question above: ropes
appear
indeed
in a
very natural and geometric way. It also answers the second question,
for
it implies that ropes are degenerations of smooth, irreducible curves,
i.e.,
it shows that ropes can be \emph{smoothed}.

\medskip

Finding out that ropes are smoothable adds to
our knowledge of the Hilbert schemes of curves. Precisely, we describe
new points which lie in the boundary of components of the Hilbert scheme that
parametrize smooth curves. This could shed
light on whether Hilbert points of smooth curves can be connected
by passing through Hilbert points of schemes without embedded points. Such is
the case for curves in $\bold P^3$.
\medskip

The way we prove the smoothing of ropes is also of independent
interest. To be precise we prove the following important fact
(see Theorem~\ref{coversmoothing} and Corollary~\ref{cover.cor}):
finite covers onto smooth
irreducible curves embedded in projective space can be deformed
to a family of $1:1$ maps.
%This was not previously
%known, to the best of the authors' knowledge.

\medskip

The relation between ropes of multiplicity $2$ and deformations of
$2:1$ morphisms had been previously studied. Fong (see ~\cite{Fong})
explored, using computational methods,
the special case of canonical ribbons and its
relation to the degeneration of canonical embeddings of smooth
curves. Recently the second author
(see ~\cite{Gonzalez03})
studied the case of ropes of multiplicity $2$ on curves of arbitrary
genus. In this context Theorem ~\ref{coversmoothing} is a more general
result as it relates ropes of arbitrary multiplicity on curves and
degenerations of embeddings to finite morphisms. In Section 3, this
result is used to prove (see Theorem  ~\ref{embsmoothing} and
Corollary ~\ref{absgensmooth}) that ropes with nonnegative arithmetic
genus and satisfying condition (*) are smoothable.

\medskip

In the second part of this article (Sections 4 and 5) we focus on the study of
ropes $\WY$ of multiplicity $3$ on $\bold
P^1$, where the general theory can be made more explicit. The conormal
bundle of  $\WY$ is a vector bundle of rank $2$, $\Cal E=\Cal O_{\bold
  P^1}(-a) \oplus \Cal O_{\bold P^1}(-b)$. Since the arithmetic genus of $\WY$ is
$a+b-2$, for a fixed arithmetic genus $p_a$ there exist ropes with different
conormal bundle. Thus, for given $p_a$, there are different families
of ropes which correspond to the nonnegative number $n=|a-b|$.
Ropes are more general when $n$ is smaller (see Proposition
~\ref{ropedegen}). This phenomenon is analogous to the one of trigonal
curves, which are stratified in the moduli of curves by the so-called
{\it Maroni invariant}.
Because of this and by analogy with the Maroni invariant of smooth
trigonal curves, we call $n$ the {\it Maroni
invariant} of $\WY$.

\medskip

In Section 4 we study the morphisms and embeddings of
ropes of multiplicity $3$ on $\bold P^1$ to projective space. The main result in this
regard is Theorem ~\ref{embedding}, which in particular tells the
smallest projective space in which all ropes with fixed conormal
bundle can be embedded. The uniform bound $N$ found depends on the
arithmetic genus and on the Maroni invariant of the rope.
We also study embeddings and morphisms from ropes to projective spaces
of dimension smaller than $N$.

\medskip

Finally, in Section 5 we apply the results of
Section 3 to smooth
ropes of multiplicity $3$ on $\bold P^1$, both as abstract schemes and
as schemes embedded in projective space, in particular, the images of the embeddings studied in Section 4 (see Theorem ~\ref{main.rope.theorem}). Knowledge of triple covers of $\bold P^1$ makes possible  to exactly characterize those ropes which satisfy (*).
These are ropes whose Maroni invariant is not too large compared to
the genus.  However, we are able to prove that not only the ropes that
satisfy (*) are smoothable but so also all other ropes with nonnegative
arithmetic genus. The reason for
this is that ropes with higher Maroni invariant can be deformed to
ropes with lower Maroni invariant (see Proposition  ~\ref{ropedegen}).

\medskip

The approach presented in this paper applies to other multiple
structures such as ropes on surfaces. For instance, in ~\cite{Enriques},
we build on the methods developed here to show the smoothing of the
so-called K3 carpets on Enriques surfaces. These ideas also throw
light on the deformation of non--locally Cohen-Macaulay
multiple structures.
The ideas used here are quite novel in this context,
comprising a combination of
techniques from new deformation theory and suitable moduli arguments.

\medskip

\noindent{\bf{Convention.}}
We work over an algebraically closed field $\bold
k$ of characteristic $0$.

\medskip
\noindent{\bf Acknowledgements:} We thank N. Mohan Kumar
for some valuable discussions and comments. We also thank E. Sernesi for
some enlightening comments.

\section{Deformation of finite covers}
\label{defcovers}

In this section we show that finite covers $X \overset{\pi} \longrightarrow Y$ of curves can always be deformed to projective embeddings. More precisely, we see that $\pi$ can be realized as a morphism $\varphi$ from $X$ to a projective space $\bold P^N$, finite onto its image, which can be obtained as a degeneration of a family of embeddings of smooth curves into $\bold P^N$.
This is essentially the content of Theorem~\ref{coversmoothing} and Corollary~\ref{cover.cor}. The proof of Theorem~\ref{coversmoothing} builds upon the ideas of the proof of \cite[Theorem 5.1]{Gonzalez03}. Theorem~\ref{coversmoothing}, however, is more general. The proof given here is also simpler and more transparent.

\begin{theorem}\label{coversmoothing}
Let $X  \overset{\pi} \longrightarrow Y$ be a cover of degree $m \geq 2$ between a smooth, irreducible, projective curve $X$ and a smooth, irreducible curve $Y$, embedded in $\bold P^r$, $r \geq 3$. Let $\Cal E$ be the trace zero module of $\pi$, i.e., let $\pi_*\Cal O_X=\Cal O_Y \oplus \Cal E$. Let $\varphi$ denote the composition of $\pi$ with the inclusion of $Y$ in $\bold P^r$.
Assume
 that
 \begin{equation*}
 h^1(\SO_Y(1))=0, h^1(\SE \otimes \SO_Y(1))=0 \ \text{and} \
 h^0(\varphi^*\Cal O_Y(1)) \geq r+1.
 \end{equation*}
 Let $\Delta=\text{Spec}\ {\bold k[\epsilon]}/{\epsilon^2}$. Then for
 every first--order infinitesimal deformation
\begin{equation*}
\WX \overset{\wphi}\longrightarrow \PP_{\Delta}^r
\end{equation*}
of $X \overset{\varphi}\longrightarrow \PP^r$, there exists a smooth
irreducible family $\SX$, proper and flat over a smooth pointed affine
curve $(T, 0)$, and a $T$--morphism $\SX \overset{\Phi}\longrightarrow
\PP_T^r$ with the following features:
\begin{enumerate}
\item the general fiber $\SX_t \overset{\Phi_t}\longrightarrow \PP^r,
  \, t \in T-0,$ is a closed immersion of a smooth irreducible
  projective curve $\SX_t$; and
\item the fiber of $\SX \overset{\Phi}\longrightarrow \PP_T^r$ over
  the tangent vector at $0 \in T$ is $\WX
  \overset{\wphi}\longrightarrow \PP_{\Delta}^r$;
in particular, the central fiber $\SX_0
\overset{\Phi_0}\longrightarrow \PP^r$ is $X
\overset{\varphi}\longrightarrow \PP^r$.
\end{enumerate}
\end{theorem}
\begin{Proof}
Let us denote $d =- \,\mathrm{deg}\, (\wedge^{m-1} \SE)$ and let $g$
be the genus of $Y$.
Note that the genus of $X$ is $g_{{}_X}=d + mg -m+ 1$.
Let $L=\pi^*\SO_Y(1)=\varphi^*\SO_Y(1)$.
The assumptions that $H^1(\SO_Y(1))=0$ and $H^1(\SE \otimes
\SO_Y(1))=0$ imply that $H^1(L)=0$.\\
Moreover we have
\begin{equation*}
H^0(L)= H^0(\SO_Y(1))\oplus H^0(\SE \otimes \SO_Y(1)).
\end{equation*}
Let us denote
\begin{equation*}
\WL=\wphi^*\SO_{\PP_{\Delta}^r}(1).
\end{equation*}
Then the restriction of $\WL$ to $X$ is $L$ and
the $\Delta$--module $\Gamma(\WL)$ is free, of rank equal to $h^0(L)$,
and $\Gamma(\WL)\otimes k[\epsilon]/\epsilon k[\epsilon]= H^0(L)$.

\medskip

\noindent We want to construct a family $\Phi$ of morphisms over $T$ such that
its fiber over the tangent vector at $0$ is $\wphi$ and its general
fiber is an embedding.
We do this in two steps. First we construct a family $(\SX,\SL)$ and
second, we construct $\Phi$ from this family.

\medskip
\noindent
\emph{Step 1.} Construction of $(\SX,\SL)$.

\smallskip
\noindent
We want to obtain a family $(\SX,\SL)$, proper and flat over a smooth pointed affine curve $(T, 0)$, whose central fiber is $(X,L)$, whose restriction to the tangent vector to $T$ at $0$ is $(\WX,\WL)$ and whose general member $(\SX_t,\SL_t)$ consists of a smooth irreducible projective curve of genus $g_{{}_X}$ and a very ample line bundle $\SL_t$.
If $g_{{}_X} \leq 2$ then $L$ is very ample, since by hypothesis $\mathrm{deg}\,L \geq 6$.
However, if $g_{{}_X} \geq 3$, then $L$ needs not be very ample,
therefore we distinguish two non--mutually exclusive cases, namely
$g_{{}_X} \geq 3$ and $L$ very ample.

\medskip
\noindent
\emph{Case 1.1} ($g_{{}_X} \geq 3$). Since we are not assuming $L$ to
be very ample, the argument that follows will be rather technical,
using the properties of the fine part of the moduli of curves and the
result of Eisenbud and Harris~\cite[5.1]{EisenbudHarris83}.

\smallskip
\noindent
We consider $\womega=\omega_{\WX/\Delta}$ and $\WL'=\WL \otimes \womega^{\otimes n}$, where $n$ is large enough so that $L'= L \otimes \omega_X^{\otimes n}$ is very ample, non--special and the complete linear series of $L'$ defines an embedding $X \rightarrow X' \subset \PP^{r'}$ that determines a smooth point $[X']$ in the corresponding Hilbert scheme.
Let $H$ be the open, smooth and irreducible subset of this Hilbert scheme consisting of the points which correspond to smooth, irreducible, non--degenerate curves $C \subset \PP^{r'}$ of degree $d'=m \, \mathrm{deg}\,\SO_Y(1)+n(2g_{{}_X}-2)$ and genus $g_{{}_X}$.
Then $[X'] \in H$.
Since $n>>0$, for every such curve $C$, $\SO_C(1)$ is non--special, the embedding of $C$ in $\PP^{r'}$ is defined by a complete series and defines a smooth point in its Hilbert scheme.
Moreover, since $L'$ is very ample and $H^1(L')=0$, $\WL'$ is very ample relative to $\Delta$ and the embedding $X \hookrightarrow \PP^{r'}$ extends to an embedding $\WX \hookrightarrow \PP_{\Delta}^{r'}$.
So the image $\WX'$ of $\WX \hookrightarrow \PP_{\Delta}^{r'}$ is a flat family over $\Delta$ that corresponds to a tangent vector to $H$ at the Hilbert point $[X']$ of $X'$.
We can take the embedding $\WX \hookrightarrow \PP_{\Delta}^{r'}$ so that this tangent vector is nonzero.
Now, since $[X']$ is a smooth point in $H$, we can take a smooth irreducible affine curve $T$ in $H$ passing through $[X']$ with tangent direction the given tangent vector.

\smallskip
\noindent
We can take the above curve in such a way that all its points except perhaps $[X']$ are placed in the open subset $U$ of $H$ constructed in the following way: $H$ admits a surjective morphism onto $\SP_{d',g_{{}_X}}$, the coarse moduli of pairs consisting of a curve of genus $g_{{}_X}$ and a line bundle of degree $d'$ on the curve.
Denote $d_1=m \,\mathrm{deg}\,\SO_Y(1)$ and consider also $\SP_{{d_1},g_{{}_X}}$, fibered over $\SM_{g_{{}_X}}^0$, the fine part of the moduli space of curves of genus $g_{{}_X}$.
Let $\SC^{(d_1)}$ be the scheme that represents the functor of relative effective Cartier divisors of relative degree $d_1$ on the universal curve $\SC_{g_{{}_X}}^0 \to \SM_{g_{{}_X}}^0$ (see \cite[4.1]{SemBourbaki232}).
We have
\begin{equation*}
d_1 - g_{{}_X}=h^0(L)-1.
\end{equation*}
So from the assumptions $h^0(L)\geq r+1$ and $r \geq 3$ we have $d_1-g_{{}_X} \geq 3$, in particular $d_1 \geq g_{{}_X}$.
Therefore the morphism $\SC^{(d_1)} \to \SP_{{d_1},g_{{}_X}}$ is surjective.
Denote $\SC = \SC_{g_{{}_X}}^0 \times_{{}_{\SM_{g_{{}_X}}^0}} \SC^{(d_1)}$. Over $\SC$ there is a universal effective relative Cartier divisor $\SD$. Consider the line bundle $\SO_{\SC}(\SD)$ and let $\SC \overset{q}\to \SC^{(d_1)}$ be the (proper and flat) projection.
Then, by the theorem of base change and cohomology, at a point $(C, D)$ of $\SC^{(d_1)}$, consisting of a curve $C$ of genus $g_{{}_X}$ and a non--special divisor $D$ of degree $d_1$ on $C$, the fiber of the coherent sheaf $R^1q_*\SO_{\SC}(\SD)$ is isomorphic to $H^1(C, D)$ and the same is true near $(C, D)$. So there is a non--empty open set $W_1$ in $\SC^{(d_1)}$ formed by pairs consisting of a curve and a divisor whose associated line bundle is non--special.
Furthermore, if we restrict to $W_1$, then the support of the cokernel of the natural map $q^*q_*(\SO_{\SC}(\SD)) \to \SO_{\SC}(\SD)$ is outside of the inverse image on an open set $W_2 \subset W_1$.\\
So we obtain an open set $W_2$ on $\SC^{(d_1)}$, formed by pairs consisting of a curve of genus $g_{{}_X}$, and an effective divisor of degree $d_1$ whose associated line bundle is non--special and globally generated and such that its associated complete linear series has dimension $h^0(L)-1 \geq 3$.
In~\cite[5.1]{EisenbudHarris83} it is proved that on a general smooth curve the general linear series of dimension $\geq 3$ has no base points and its associated map to projective space is a closed immersion.
Moreover, $\SC^{(d_1)}$ is irreducible so $W_1$ dominates $\SM_{g_{{}_X}}^0$.
Therefore, the set $W_2$ is non--empty, and shrinking $W_2$ so that $q_*(\SO_{\SC}(\SD))$ is free of rank $h^0(L)$ on $W_2$, we have a $W_2$--morphism
$\, \SC_{g_{{}_X}}^0 \times_{{}_{\SM_{g_{{}_X}}^0}} W_2  \to \PP_{W_2}^{h^0(L)-1}$.
Now (see e.g.~\cite[4.6.7]{EGA3-1}),  the points of $W_2$ such that the morphism induced on the fiber over the point is a closed immersion form an open set $W$ in $W_2$.
So we obtain an open set $W$ in $\SC^{(d_1)}$ formed by pairs consisting of a curve of genus $g_{{}_X}$ and an effective divisor of degree $d_1$ whose associated line bundle is non--special and very ample.
Also by~\cite[5.1]{EisenbudHarris83}, as $h^0(L)-1  \geq 3$, we see that the open set $W$ is non--empty.
Now, since $\SC^{(d_1)}$ is irreducible and $\SC^{(d_1)} \to
\SP_{{d_1},g_{{}_X}}$ is surjective, we also obtain a non--empty open
set in $\SP_{d_1,g_{{}_X}}$ formed by pairs consisting of a curve and
a very ample non--special line bundle with as many global sections as
$L$.
Moreover, twisting by $\omega^{\otimes n}$ we have an isomorphism between $\SP_{d_1,g_{{}_X}}$ and $\SP_{d',g_{{}_X}}$.
So we take the open set $\,U \subset H$, which is the  inverse image of the considered open set in $\SP_{d',g_{{}_X}}$.\\
Let $0 \in T$ denote the point corresponding to $X'$.
Now, we take the curve $T$ so that $T-\{0\} \subset U$.
Over the pointed affine curve $(T, 0)$ we have a proper flat polarized family $(\SX,\SL')$ whose fibers  over $0$ and over the tangent vector to $T$ at $0$ are $(X, L')$ and $(\WX,\WL')$ respectively.\\
Now, twisting by $\omega_{{}_{\SX/T}}^{-n}$, we obtain a family $(\SX,\SL)$, proper and flat over $T$, whose central fiber is $(X,L)$, whose restriction to the tangent vector to $T$ at $0$ is $(\WX,\WL)$ and whose general member $(\SX_t,\SL_t)$ consists of a smooth irreducible projective curve of genus $g_{{}_X}$ and a very ample line bundle $\SL_t$. Furthermore, the total family $\SX$ is smooth and irreducible and $\SL_t$ is non--special and has as many global sections as $L$ and degree $d_1=\,\mathrm{deg}\,L$.

\medskip
\noindent
\emph{Case 1.2} ($L$ very ample). This happens when $g_{{}_X} \leq 2$
or in most cases when $g_{{}_X}\geq 3$. In this case the argument
above can be drastically simplified. Being $L$ very ample and non--special,
 its complete linear series defines an embedding $X
\rightarrow X' \subset \PP^{d_1-g_{{}_X}}$ that determines a smooth
point $[X']$ in the corresponding Hilbert scheme. Now in the open,
smooth subset $H$  of this Hilbert scheme consisting of the points which correspond to smooth, irreducible, non--degenerate curves $C \subset \PP^{d_1-g_{{}_X}}$ of degree $d_1=\mathrm{deg}\,L$ and genus $g_{{}_X}$ with non--special hyperplane section, we have a point $[X']$ and a nonzero tangent vector corresponding to an embedding $\WX \hookrightarrow \PP_{\Delta}^{d_1-g_{{}_X}}$, extending the embedding $X \hookrightarrow \PP^{d_1-g_{{}_X}}$, given by the very ample line bundle $\WL$.
So we can take a smooth irreducible affine curve $T$ in $H$ passing
through $[X']$ and such that its tangent direction is the given
tangent vector.
Then the pullback to $T$ of the universal family provides
a family $(\SX,\SL)$, proper and flat over
$T$, whose central fiber is $(X,L)$, whose restriction to the tangent
vector to $T$ at $0$ is $(\WX,\WL)$ and whose general member
$(\SX_t,\SL_t)$ consists of a smooth irreducible projective curve of
genus $g_{{}_X}$ and a very ample, non--special line bundle $\SL_t$
with as many global sections as $L$ and degree $d_1=\mathrm{ deg }L$.
%Observe that since $L$ is very ample we do not need to
%use~\cite[5.1]{EisenbudHarris83}, because it is well known that a
%general subseries of dimension $\geq 3$, corresponding to a very
%ample line bundle, is very ample.

\medskip
\noindent
\emph{Step 2.} Construction of $\Phi$.

\smallskip
\noindent
We will construct, after shrinking $T$ if necessary,  a $T$--morphism
\begin{equation*}
\SX \overset{\Phi}\to \PP_T^r
\end{equation*}
whose fiber over the tangent vector to $T$ at $0$ is the morphism $\WX \overset{\wphi}\to \PP_{\Delta}^r$ in the statement and whose general fiber $\SX_t \overset{\varphi_t} \to \PP^r$ for $t \neq 0$ is a closed immersion given by a $r$--dimensional linear subseries of $H^0(\SL_t)$.

\medskip

\noindent First we look at the properties of $\SL$ and of its module of global sections.  Let $\SX \overset{p}\to T$ be the (proper and flat) morphism from $\SX$ to $T$.
The facts that $p$ is proper, $\SL$ is flat over $T$, and $h^0(\SL_t)=h^0(L)$ and $H^1(\SX_t,\SL_t)=0$ for every $t \in T$, imply that $p_* \SL$ is a locally free sheaf of rank $h^0(L)$ on $T=\mathrm{Spec}\, R$ and that ``the formation of $p_*$ commutes with base extension''. So we have $\Gamma(\SL) \otimes_R k[\epsilon]/\epsilon k[\epsilon] = \Gamma(L)$, $\Gamma(\SL) \otimes_R k[\epsilon]/\epsilon^2 = \Gamma(\WL)$ and $\Gamma(\SL)\otimes_R k(t)=H^0(\SX_t,\SL_t)$ for every point $t \in T$, $t \neq 0$.
After shrinking $T$, we can assume that $\Gamma(\SL)$ is a free $R$--module of rank $h^0(L)=n+1$ that induces a $T$--morphism
\begin{equation*}\label{Psi}
\SX \xrightarrow{\Psi} \PP_{T}^n.
\end{equation*}

\medskip

\noindent Let $\WPsi$ the restriction of $\Psi$ to $\Delta$. Since $\wphi$ is induced by a set of global sections of $\Gamma(\WL)$ and $\Gamma(\SL) \otimes_R k[\epsilon] = \Gamma(\WL)$, we have that there exists a linear projection
\begin{equation*}\label{wrho}
\PP_{\Delta}^n \overset{\wrho}\dashrightarrow \PP_{\Delta}^r,
\end{equation*}
such that $\wphi = \wrho \circ \WPsi$.
Now let $t_1 \in T$, $t_1 \neq 0$. Since $\Gamma(\SL)\otimes_R k(t_1)=H^0(\SX_{t_1},\SL_{t_1})$, the restriction $\Psi_{t_1}$ of $\Psi$ to $t_1$ is induced by $H^0(\SL_{t_1})$. By Step 1, $\SL_{t_1}$ is very ample, so $\Psi_{t_1}$ is a closed immersion into $\PP^n_{t_1}$. Since $\SX_{t_1}$ is a curve and $3 \leq r \leq h^0(\SL_t)-1$, it is well known that for a general projection $\rho_{t_1}$ the composition $\rho_{t_1} \circ \Psi_{t_1}$ is again a closed immersion.

\medskip

\noindent Now let $\Sigma=\Delta \cup \{t_1\}$. The pair $(\wrho, \rho_{t_1})$ is represented by a matrix with coefficients in $\SO_\Sigma$. Lifting them to $R$ we obtain a linear projection
\begin{equation*}\label{projection}
\PP_{T}^n \overset{\rho}\dashrightarrow \PP_{T}^r
\end{equation*}
which restricts to $\wrho$ and to $\rho_{t_1}$. We define $\Phi$ as the composition $\rho \circ \Psi$. After maybe shrinking $T$ we may assume that $\Phi$ is $T$-morphism.
%We know that $\Phi_0$ and $\Phi_{t_1}$ are morphisms. In general $\Phi$ is only map, thus we have $ev: \S0_{\SX}^{r+1} \to \SL$, which is surjective at the fiber at $0$ and $t_1$. Then $ev$ is surjective in a neighborhood of $0$ and $t_1$.
After maybe shrinking $T$ again we see that $\Phi$ is a closed immersion when restricted to every $t \in T$, $t \neq 0$ (see  ~\cite[4.6.7]{EGA3-1}). By construction the restriction of $\Phi$ to $\Delta$
is $\WX \overset{\wphi}\to \PP_{\Delta}^r$.
\end{Proof}

\begin{cor}\label{cover.cor}
Let $X \overset{\pi} \longrightarrow Y$ be a finite cover between a smooth, irreducible curve $X$ and a smooth, irreducible curve $Y$.  Then $\pi$ can be realized as limit of projective embeddings of smooth curves. More precisely, there are embeddings $i$ of $Y$ in projective space so that the composition $\varphi=i \circ \pi$ is the limit of embeddings of smooth curves to projective space.
\end{cor}

\begin{proof}
Let $\SE$ be the trace zero module of $\pi$. There are always embeddings $i$ of $Y$, induced by complete linear series,  into a projective space of sufficiently high dimension  so that $H^1(\SO_Y(1))=0$, $H^1(\SE \otimes \SO_Y(1))=0$. Then we can apply Theorem~\ref{coversmoothing}.
\end{proof}

\section{Smoothing of ropes on curves}
\label{generalsmooth}

In this section we show that ropes of arbitrary multiplicity supported on curves of arbitrary genus can be smoothed, i.e., can be obtained as flat limits of smooth curves, provided they satisfy certain mild conditions which are quite natural from a geometric point of view. We start recalling the definition of ropes:

\begin{definition}\label{defrope}
Let $Y$ be a reduced connected scheme and let $\SE$ be a locally free sheaf of rank $m-1$ on $Y$.
A rope of multiplicity $m$ or an $m$-rope, for short, on $Y$ with conormal bundle $\SE$ is a scheme $\WY$ with ${\WY}_{\mathrm{red}}=Y,$ such that
\begin{enumerate}
\item
$\SI_{Y, \WY}^2=0$ and
\item
$\SI_{Y, \WY} \simeq \SE$ as $\SO_{Y}$--modules.
\end{enumerate}
When $\SE$ is a line bundle, $\WY$ is called a ribbon on $Y$.
\end{definition}

By smoothing a rope $\WY$ we mean finding a flat, integral  family $\SY$ of schemes over a smooth affine curve $T$, such that over a point $0 \in T$, $\SY_0=\WY$ and over the remaining points $t$ of $T$, $\SY_t$ is a smooth, irreducible variety. We will prove that a rope $\WY$ of multiplicity $m$ on $Y$ is smoothable if its conormal bundle $\Cal E$ is the trace zero module of a smooth, irreducible cover $\pi$ of $Y$ of degree $m$. We will show then that $\WY$ appears as the limit of the images of a family of projective embeddings degenerating to $\pi$.
Thus,  we will first consider $\WY$ embedded in projective space.  Then we use the fact that $\WY$ arises as the central fiber of the image of a first--order infinitesimal deformation of the composition of $\pi$ with the inclusion of $Y$ in projective space $\bold P^N$.  Before doing all these, we need  to recall the spaces parametrizing ropes and morphisms from ropes of fixed conormal bundle $\SE$ to projective space:

\begin{theorem}\label{paropes}
Let $Y$ be a reduced connected scheme and let $\Cal E$ be a vector bundle of rank $m-1$ on $Y$.

\begin{enumerate}
\item
A rope $\WY$ on $Y$ with conormal bundle $\SE$ is determined by
an element $[e_{\WY}] \in \mathrm{Ext}_Y^1(\Omega_Y, \SE)$. The rope
$\WY$ is split (i.e., $Y \hookrightarrow \WY$ admits a retraction) if
and only if $[e_{\WY}]=0$.

\smallskip

\noindent Assume furthermore that $Y$ is a smooth variety  and let $Y \overset{i}\hookrightarrow Z$ be a closed immersion in another smooth variety $Z$.

\smallskip

\item There is a one--to--one correspondence between pairs $(\WY, \wi)$, where $\WY$ is a rope on $Y$ with conormal bundle $\,\SE$ and $\, \WY \overset{\wi}{\to} Z$ is a morphism extending $Y \overset{i}{\hookrightarrow} Z$, and elements $\tau \in \mathrm{Hom}(\SN_{Y,Z}^*,\SE)$.

\item  If $\tau$ corresponds to a pair $(\WY, \wi)$, $\wi$ is an embedding if and only if  $\tau$ is a surjective homomorphism.

\item If $\tau$ corresponds to a pair $(\WY, \wi)$, then $\tau$ is mapped by the connecting homomorphism $\delta$
onto $[e_{\WY}]$.

\end{enumerate}
\end{theorem}

\begin{proof}
For proofs of the statements see ~\cite{BE} and \cite [Proposition 2.1]{Gonzalez03}.
\end{proof}

The construction of the first infinitesimal deformation mentioned
above and the relation of ropes and  infinitesimal deformations of
finite morphisms was done in ~\cite[Theorem 3.10]{Gonzalez03}. We
state now a direct consequence of this theorem:

\begin{theorem}\label{infsmoothing}  Let $\WY \subset \bold P^r$ be a rope of multiplicity $m$ on a smooth irreducible curve $Y$  and let  $\SE$ be the  conormal bundle of $\WY$.
Assume that there exists a smooth irreducible cover $X \overset{\pi}\longrightarrow Y$ of degree $m$ such that $\pi_*\SO_X=\SO_{Y}\oplus \SE$.
Let $X \overset{\varphi}\longrightarrow \PP^r$ be the morphism obtained by composing $\pi$ with the inclusion of $\,Y$ in $\PP^r$.
Then $\WY$ is the central fiber of the image of some first--order infinitesimal deformation of $\varphi$.
\end{theorem}

Now we use Theorems~\ref{coversmoothing} and ~\ref{infsmoothing}
 to show that $\WY$ is the limit of the images of a family of embeddings $\Phi_t$ of smooth curves, degenerating to $\varphi$. Precisely, we want to extend the infinitesimal deformation of $\varphi$ in such a way that, if we call the image of the family of morphisms $\SY \subset \PP^N \times T$, then $\SY_0=\WY$. All this is done in the next theorem:

\begin{theorem}\label{embsmoothing}
Let $\WY$ be a rope of multiplicity $m$ and nonnegative arithmetic genus, embedded in $\PP^r$
($r \geq 3$) as non--degenerate subscheme, and supported on a smooth irreducible projective curve $Y$. Let $\SE$ be the conormal bundle of $\WY$.  Assume that $H^1(\SO_Y(1))=0$ and $H^1(\SE \otimes \SO_Y(1))=0$.

\smallskip

\noindent If
there exists a smooth irreducible cover $X \overset{\pi}\longrightarrow Y$ such that $\pi_*\SO_X=\SO_{Y}\oplus \SE$, then there exists a family of morphisms $\Phi$ over an affine curve $T$ as described in Theorem~\ref{coversmoothing} such that the image $\Cal Y$ of $\Phi$  is a closed integral subscheme $\SY \subset \PP_T^r$, flat over $T$, with the following features:
\begin{enumerate}
\item the general fiber $\,\SY_t , \, t \in T-0,$ is a smooth irreducible projective non--degenerate curve with non--special hyperplane section in $\PP^r$,
\item the central fiber $\,\SY_0 \subset \PP^r$ is $\,\WY \subset \PP^r$.
\end{enumerate}
%\end{theorem}

\begin{remark} {\rm The hypothesis in Theorem ~\ref{embsmoothing}
that $\WY$ be embedded in $\PP^r$
as non--degenerate subscheme can be relaxed. Indeed, one
only needs to ask that $h^0(\pi^*\SO_Y(1)) \geq r+1$ (see
~\eqref{h0L}).}
\end{remark}

\begin{proof}(of Theorem ~\ref{embsmoothing})
We use the notations of the proof of the Theorem~\ref{coversmoothing}.
For the line bundle $L=\pi^*\SO_Y(1)$ on $X$ we have
\begin{equation*}
H^0(L)= H^0(\SO_Y(1))\oplus H^0(\SE \otimes \SO_Y(1)).
\end{equation*}
Moreover, from the sequence
\begin{equation*}\label{rope-algebra}
\xymatrix@1{
0 \ar[r] & \SE  \ar[r] & \SO_{\WY}  \ar[r] & \SO_Y \ar[r] & 0,}
\end{equation*}
twisted by $\SO_{\WY}(1)$ and the assumption $H^1(\SE \otimes \SO_Y(1))=0$ we obtain an isomorphism
\begin{equation*}\label{h0-rope-algebra}
H^0(\SO_{\WY}(1))=H^0(\SO_Y(1))\oplus H^0(\SE \otimes \SO_Y(1)).
\end{equation*}
By assumption $\WY$ is embedded in $\PP^r$ as a non--degenerate subscheme, so $H^0(\SO_{\PP^r}(1)) \subset H^0(\SO_{\WY}(1))$.\\
Thus  we obtain the inequality
\begin{equation}\label{h0L}
h^0(L) \geq r+1.
\end{equation}
Moreover from Theorem~\ref{infsmoothing} we know that there exists a first--order infinitesimal deformation
\begin{equation*}
\WX \overset{\wphi}\to \PP_{\Delta}^r
\end{equation*}
of $\varphi$ such that the central fiber of the image of $\wphi$ is equal to the rope $\WY$.\\
Therefore
there exist a family $\SX \to T$ and a $T$--morphism $\SX \overset{\Phi}\to \PP_T^r$ as
in Theorem~\ref{coversmoothing}.
Let $\SY$ be the image of the $T$--morphism $\SX \overset{\Phi}\to \PP_T^r$.
The total family $\SX$ is smooth and irreducible so $\SY$ is integral.
Furthermore, $\Phi$ is a closed immersion over $T-0$ since, by Theorem~\ref{coversmoothing},  $\Phi_t$ is a closed immersion for every $t \in T-0$ (see e.g. \cite[4.6.7]{EGA3-1}).
Therefore for $t \in T-0$ we have the equality $\SY_t = \,\mathrm{im}\,(\Phi_t)$. Since $\SX_t$ is smooth, this proves (1).
Finally, the fact that $T$ is an integral smooth curve and $\SY$ is integral and dominates $T$ imply that $\SY$ is flat over $T$.
So the fiber $\SY_0$ of $\SY$ at $0 \in T$ is the flat limit of the images of $\SX_t \overset{\Phi_t} \to \PP^r$ for $t \neq 0$.
Moreover, this fiber $\SY_0$ contains the central fiber $\,(\mathrm{im}\,\widetilde{\varphi})_0$ of the image of $\wphi$.
Since $\WY$ has conormal bundle $\SE$ and $X$ has trace zero module $\SE$, the genus of $X$ and the arithmetic genus of $\WY$ are equal.
Then $\SY_0$ and $(\mathrm{im}\,\widetilde{\varphi})_0$ have the same degree and the same arithmetic genus, so they are equal.
\end{proof}

\end{theorem}

\begin{remark}\label{manysmooth}
{\rm In the vast majority of  cases,
for a fixed inclusion $\WY \overset{\tilde i} \hookrightarrow  \bold P^r$,  there are many different smoothings of $\WY$:
\begin{enumerate}
\item First, there will be many possible covers $\pi$ to choose.
\item  Second, once $\pi$ is chosen and $\varphi$ is therefore fixed, the element $\tau \in \mathrm{Hom}(\SN_{Y,Z}^*,\SE)$ corresponding to $(\tilde Y, \tilde i)$ may have, in most cases, many different liftings to $H^0(\Cal N_\varphi)$. Liftings $\mu$ of $\tau$ to $H^0(\Cal N_\varphi)$ correspond to first--order infinitesimal deformations $\widetilde \varphi$ of $\varphi$, such that (im$\widetilde \varphi)_0=\WY$. For details on this, see ~\cite[Section 3]{Gonzalez03}.
\item Third, there are many different ways of extending $\widetilde \varphi$ to a family of morphisms $\Phi$, as the proof of Theorem~\ref{coversmoothing} shows.
\end{enumerate} }
\end{remark}

\begin{cor}\label{absgensmooth}
Let $\WY$ be a rope of multiplicity $m$ on a smooth irreducible projective curve $Y$,  with conormal bundle $\SE$ and nonnegative arithmetic genus. If there is a smooth connected cover $X \overset{\pi}\longrightarrow Y$ such that $\pi_*\SO_X=\SO_{Y}\oplus \SE$, then $\WY$ is smoothable.
\end{cor}

\begin{proof}
Since $\WY$ is a proper scheme of dimension $1$, it is projective. Taking $r$ sufficiently large we can embed $\WY$ as a non--degenerate subscheme of $\bold P^r$ and at the same time, make $H^1(\Cal O_Y(1))=H^1(\SE\otimes \SO_Y(1))=0$. Then the corollary follows from Theorem~\ref{embsmoothing}.
\end{proof}

\section{Embeddings of multiplicity $3$ ropes in projective space}
\label{embeddings}

Proper schemes of dimension $1$ are projective, so in particular, any multiplicity $3$ rope $\WY$ on $\bold P^1$ can be embedded in a sufficiently large projective space $\bold P^N$. The main purpose of this section is to find, given $\WY$, the smallest $N$ so that $\WY$ can be embedded in $\bold P^N$.

%
%\smallskip
%
The conormal bundle of a rope $\WY$ of multiplicity $3$  is a vector
bundle $\Cal E$ of rank $2$.  The bundle $\Cal E$ splits as a direct sum $\Cal O_{\bold P^1}(-a) \oplus \Cal O_{\bold P^1}(-b)$, where $a$ and $b$ are integers and $a \geq b$. The arithmetic genus of $\WY$ is $a+b-2$.
Thus, if we set an arithmetic genus $g$, there are ropes of genus $g$ with different conormal bundles. We can set a hierarchy among them:

\begin{defi}\label{balanced}
Let $\WY$ be a rope on $\bold P^1$ with conormal bundle $\Cal E=\Cal O_{\bold P^1}(-a) \oplus \Cal O_{\bold P^1}(-b)$, $a \geq b$.
\begin{enumerate}
\item By analogy with smooth trigonal curves, we say that the Maroni invariant of $\WY$ is $n={a-b}$

\item
We say that $\WY$ is balanced if $n=0$ or $n=1$. Let $\WY_1$ and $\WY_2$ have
the same arithmetic genus and let $n_i$  be the Maroni invariant of $\WY_i$. We say that $\WY_1$ is more balanced than $\WY_2$ if $n_1 \leq n_2$.
\end{enumerate}
\end{defi}

In Section~\ref{triplesmooth} (see Proposition~\ref{ropedegen}), we will show that a less balanced rope is a degeneration of more balanced ropes.

\medskip

Now we proceed to describe the embeddings of ropes of multiplicity $3$ on $\bold P^1$.
We start looking for embeddings supported on a rational normal curve
lying in the same projective space where the rope is
embedded. Obviously, a rope of multiplicity $3$ cannot be embedded in
$\bold P^2$. There is not much room in $\bold P^3$ either. Indeed,
the only $3$-ropes in $\bold P^3$, supported on the twisted cubic are
the non--split ropes on $\bold P^1$ with conormal bundle $\Cal O_{\bold P^1}(-5)
\oplus \Cal O_{\bold P^1}(-5)$. This is because the conormal bundle of
the twisted cubic is $\Cal O_{\bold P^1}(-5) \oplus \Cal O_{\bold
  P^1}(-5)$ and the map $\delta$ in ~\eqref{embseq} is, in this case, an
isomorphism.  Thus we start looking for embeddings in projective spaces of dimension $4$ or higher. More precisely, we will
find the smallest $N$ so that all ropes with fixed conormal bundle $\mathcal E$ can be embedded in  $\bold P^N$.

\begin{theorem}\label{embedding}
Let $\WY$ be a rope on $Y=\bold P^1$ of arithmetic genus $g$, Maroni
invariant $n$ and conormal bundle $\Cal E=\Cal O_{\bold P^1}(-a)
\oplus \Cal O_{\bold P^1}(-b)$, $a \geq b$. Let $N_0=\max\{4,a-1\}$
(i.e., $N_0=\max\{4, \frac {g+n}{2}\}$).
Let $Y  \overset{i} \hookrightarrow \bold P^N$ be an embedding of $Y$ as a rational normal curve in $\bold P^N$. If $N \geq N_0$,
then there are embeddings
$\WY  \overset{\tilde i} \hookrightarrow  \PP^N$ that extend $i$. In particular, $\WY$ can be embedded as a subscheme of $\bold P^N$, supported on a rational normal curve of degree $N$.
\begin{remark}\label{linearly.normal}
{\rm Under the hypotheses of Theorem~\ref{embedding}
\begin{enumerate}
\item $H^1(\SE \otimes \SO_{\PP^1}(N))=0$.
\item $H^0(\SE \otimes \SO_{\PP^1}(N))=0$ if and only if
$\SE=\SO_{\PP^1}(-a) \oplus \SO_{\PP^1}(-a)$ (i.e., if $g$ is even and
$n=0$) and $N=a-1$ (i.e., $N=\frac g 2$).
\end{enumerate}

In particular, $\tilde i$ embeds $\WY$ as a linearly normal subscheme of $\bold P^N$ if and only if
$\SE=\SO_{\PP^1}(-a) \oplus \SO_{\PP^1}(-a)$ and $N=a-1$.}
\end{remark}
\end{theorem}

In order to prove Theorem~\ref{embedding} we need the following

\begin{lemma}\label{split.emb}
Let $Y$ be an irreducible projective variety and let $\SE$ and $\SF$ be locally free sheaves of finite rank on $Y$.
\begin{enumerate}
\item The surjective homomorphisms from $\SF$ to $\SE$ form an open set of
$\, \mathrm{Hom}(\SF,\SE)$ which is the complement of an algebraic cone.
\item Consider an extension of a coherent sheaf $\SF'$ by $\SF$.
Let $\, \mathrm{Hom}(\SF,\SE) \overset{\delta}\to \mathrm{Ext}^1(\SF',\SE)$ be the connecting map induced by the extension.
If the class of the split extension lifts to an epimorphism, then every class in the image of $\,\delta$ can be lifted to an epimorphism.
\end{enumerate}
In particular, if $Y \overset{i} \hookrightarrow Z$ is an inclusion of smooth projective varieties, and the split rope on $Y$ with conormal bundle $\Cal E$ admits an embedding into $Z$ extending $i$, then every rope on $Y$ with conormal bundle $\Cal E$ which admits a morphism to $Z$ extending $i$  also admits an embedding into $Z$ extending $i$.
\end{lemma}

\begin{proof}
The lemma is a generalization of ~\cite[Lemma 4.1]{Gonzalez03}  and its proof follows almost word by word the proof given there. The conclusion about the embeddings of ropes follows from Part (2) and Theorem~\ref{paropes}.
\end{proof}

\begin{Proof} (of Theorem~\ref{embedding})
Recall the sequence
\begin{equation}\label{embseq}
 \text{Hom}(\Omega_{\bold
P^{N}}|_{\bold P^1},\SE) \overset{\gamma} \longrightarrow \text{Hom}(\mathcal I/\mathcal I^2,
\mathcal E) \overset{\delta} \longrightarrow
\text{Ext}^1(\Omega_Y,\mathcal E) \longrightarrow  \text{Ext}^1(\Omega_{\bold
P^{N}}|_{\bold P^1},\SE).
\end{equation}
The bundle $\Omega_{\bold
P^{N}}|_{\bold P^1}=\mathcal O_{\bold P^1}(-N-1)^{\oplus N}$, hence $$\text{Ext}^1(\Omega_{\bold
P^{N}}|_{\bold P^1},\SE)=
H^1(\Cal O_{\bold P^1}(N+1-a)^{\oplus N} \oplus \Cal O_{\bold P^1}(N+1-b)^{\oplus N})=0.$$
This means that  any $\WY$ with conormal bundle $\SE$ admits a morphism to $\bold P^N$ extending $i$.
We are going to prove that, among these morphisms, one can actually find embeddings.
By Lemma~\ref{split.emb},
it is enough to show that the split rope with conormal $\SE$ possesses such an embedding, i.e., it is enough to show that there exists an epimorphism $\tau \in \text{Hom}(\mathcal I/\mathcal I^2,
\mathcal E)$ such that $\delta(\tau)=0$. This is the same as finding an
epimorphism $\tau$ in the image of $\gamma$.

\smallskip

Let $\mathcal I$ denote the ideal sheaf of the image of $Y$ in $\bold P^N$.
Then the conormal bundle $$\mathcal I/\mathcal I^2=\mathcal O_{\bold
  P^1}(-N-2)^{\oplus N-1}.$$
Moreover
$\mathcal I/\mathcal I^2$
fits in the
following exact sequence:

$$0 \longrightarrow \mathcal I/\mathcal I^2 \longrightarrow \Omega_{\bold P^{N}}|_Y
\longrightarrow \Omega_Y \longrightarrow 0 \ ,$$

which is the same as

$$0 \longrightarrow \mathcal O_{\bold P^1}(-N-2)^{\oplus N-1} \overset{\alpha} \longrightarrow
\mathcal O_{\bold P^1}(-N-1)^{\oplus N} \overset{\beta} \longrightarrow \mathcal O_{\bold P^1}(-2)
\longrightarrow 0 \ .
$$

If we fix a suitable ``basis'' for each bundle in $(*)$,
the
matrix of $\beta$ is $$\begin{pmatrix} X_0^{N-1} & X_0^{N-2}X_1 & \cdots &
X_0X_1^{N-2} & X_1^{N-1} \end{pmatrix} $$
and the matrix of $\alpha$ is the $N \times (N-1)$ matrix
$$\begin{pmatrix}  X_1 & 0 & 0 & 0 & \cdots &0 &  0 \cr
   -X_0 & X_1 & 0 & 0 &\cdots &0 &  0 \cr
    0 & -X_0 & X_1 & 0 & \cdots &0 &  0 \cr
    \cdots & \cdots & \cdots & \cdots & \cdots & \cdots & \cdots \cr
    0 & 0 & 0 & 0 & \cdots & -X_0 & X_1 \cr
    0 & 0 & 0 & 0 & \cdots & 0 & -X_0 \end{pmatrix} ,$$
 \noindent where $X_0, X_1$ are homogeneous coordinates of $\bold P^1$.
On the other hand,  an arbitrary element of $\text{Hom}(\Omega_{\bold
  P^{N}}|_Y, \mathcal E)$ has a matrix of the form
$$\begin{pmatrix} a_1 & a_2 & \cdots & a_{N} \cr
b_1 & b_2 & \cdots & b_{N} \end{pmatrix},$$
\noindent where the $a_i$s and $b_j$s are homogeneous forms on $\bold P^1$.
The forms $a_1, \cdots, a_{N}$ have degree $N-a+1$ (which is a nonnegative integer) and $b_1, \cdots, b_{N}$ have degree $N-b+1$.

Then the image of $\gamma$, which is the kernel of $\delta$,
consists of homomorphisms in
$\text{Hom}(\mathcal I/\mathcal I^2,
\mathcal E)$ whose matrices  have size $2 \times (N-1)$ and are of the form

$$M=\begin{pmatrix} a_1X_1-a_2X_0 & a_2X_1-a_3X_0 & \cdots &
a_{N-1}X_1-a_{N}X_0 \cr
b_1X_1-b_2X_0 & b_2X_1-b_3X_0 & \cdots &
b_{N-1}X_1-b_{N}X_0 \end{pmatrix} \ ,$$
where $a_1, \cdots, a_{N}, b_1, \cdots, b_{N}$ are as above.
Then, as far as $N \geq 4$, we can always find
$a_i$s and $b_j$s so that $M$ has maximal
rank $2$
at every point of $\bold P^1$ and corresponds therefore to an element of
$\text{Hom}(\mathcal I/\mathcal I^2,
\mathcal E)$ which is an epimorphism and whose image by $\delta$ is $0$.
\end{Proof}\\

In the next theorem we show that the bound on $N$ obtained in Theorem~\ref{embedding} is the best possible that allows all ropes $\WY$ with a fixed conormal bundle $\SE$ to be embedded in $\bold P^N$. Next theorem also characterizes, among other things, what ropes can be embedded in smaller projective spaces:

\begin{theorem}\label{lower.embedding.even}
Let $\WY$ be a rope on $Y=\bold P^1$ with conormal bundle $\Cal E=\Cal O_{\bold P^1}(-a) \oplus \Cal O_{\bold P^1}(-b)$, $a \geq b$.
Let $Y  \overset{i} \hookrightarrow \bold P^N$ be an embedding of $Y$ as a rational normal curve in $\bold P^N$.
\begin{enumerate}
\item
 If $N =a-2$ and $\WY$ does not admit a projection onto a ribbon  of conormal bundle $\SO_{\bold P^1}(-b)$, then there exist embeddings $\widetilde Y \overset{\tilde i}  \hookrightarrow  \bold P^N$ which extend  $i$.
\item
If $N =a-2$ and $\WY$ admits a projection onto a ribbon $\widehat Y$ of conormal $\SO_{\bold P^1}(-b)$, then there exist morphisms $\widetilde Y \overset{\tilde i}  \longrightarrow  \bold P^N$ extending  $i$ and all of them factor through a projection onto a ribbon whose conormal bundle is a subsheaf of $\SO_{\bold P^1}(-b)$.
 \item
If $N \leq a-3$, then there exist morphisms
$\widetilde Y \overset{\tilde i}  \longrightarrow  \bold P^N$ extending $i$ only if $\widetilde Y$ admits a projection $p$ onto a ribbon $\widehat Y$ whose conormal bundle is a subsheaf of $\SO_{\bold P^1}(-b)$. These morphisms factor through $\widehat Y$.

\end{enumerate}

In particular, if $N < a-1$, the rope $\WY$ can be embedded in $\bold P^N$ supported on a rational normal curve of degree $N$ if and only if $N=a-2$ and $\widetilde Y$ does not admit a projection onto a ribbon with conormal bundle $\SO_{\bold P^1}(-b)$.
\end{theorem}

\begin{proof}
Let $\Cal I$ be ideal sheaf of $i(Y)$ in $\bold P^N$. Recall that the
conormal bundle of $i(Y)$ in $\bold P^N$ is $\Cal I/\Cal I^2=\Cal
O_{\bold P^1}(-N-2)^{\oplus N-1}$.
We see first under what conditions $i$ can be extended to a morphism $\widetilde Y \overset{ \tilde i } \longrightarrow \bold P^N$,  for some rope $\widetilde Y$  with conormal bundle $\Cal E$.  Recall that, by Theorem ~\ref{paropes}, the space Hom$(\Cal I/\Cal I^2, \Cal E)$ parametrizes pairs $(\widetilde Y, \tilde i)$, where $\tilde i$ extends $i$.

\smallskip

If $N \leq b-3$, then Hom$(\Cal I/\Cal I^2, \Cal E)=0$,
so $i$ can be extended if and only if $\widetilde Y$ is the split rope. In this case it is clear that $\tilde i$ factorizes through the retraction of the split rope onto $Y$.
If $b-2 \leq N \leq a-3$, any homomorphism from $\Cal I/\Cal I^2$ to $\Cal E$ is the composition of a homomorphism from $\Cal I/\Cal I^2$ to $\Cal O_{\bold P^1}(-b)$, followed by the bundle embedding
$\Cal O_{\bold P^1}(-b) \hookrightarrow \Cal E$. This means that any extension of $i$ to a rope $\WY$ with conormal bundle $\SE$ factors through an extension (actually, an embedding, provided the corresponding homomorphism from $\Cal I/\Cal I^2$ to $\SE$ is nonzero)  of $i$ to a ribbon $\widehat Y$, whose conormal bundle is a subsheaf of $\Cal O_{\bold P^1}(-b)$. In particular, $\widetilde Y$ admits a projection onto $\widehat Y$. This proves 3.

\smallskip

We prove now 1. If $N=a-2$, then Ext$^1(\Omega_{\bold P^N}|_{\bold
  P^1}, \SE)=0$, therefore, given any rope $\widetilde Y$ with
conormal bundle $\Cal E$, the morphism $i$ can be extended to
$\widetilde Y$. Let us now fix $\widetilde Y$, which corresponds to
some $\zeta \in \text{Ext}^1(\Omega_{\bold P^1}, \SE)$, and look at the
possible extensions of $i$ to $\widetilde Y$, which correspond to
elements $\tau \in \text{Hom}(\Cal I/\Cal I^2, \Cal E)$ mapping to
$\zeta$.  As argued in the proof of Theorem~\ref{embedding}, $\tau$ can
be thought as a $2 \times (N-1)$ matrix, in which the entries of the
first row are constant and the entries of the second row are
homogeneous forms of degree $a-b$. If all the entries of the first row
are $0$, then $\tau$ corresponds to an extension $\widetilde Y
\overset{\tilde i} \longrightarrow \bold P^N$ that factors through a
ribbon whose conormal bundle is  a subsheaf of $\Cal O_{\bold
  P^1}(-b)$.  On the other hand, an element in Hom$(\Omega_{\bold
  P^N}|_{\bold P^1}, \SE)$ corresponds to a $2 \times N$ matrix with
zero entries in its first row and homogeneous forms of degree $a-b-1$
  in its second row. Arguing as in the proof of
  Theorem~\ref{embedding}, we can then add to $\tau$ a suitable element
  of the kernel of $\delta$, and obtain an element $\tau'$ in
  $\text{Hom}(\Cal I/\Cal I^2, \Cal E)$, such that
  $\delta(\tau')=\zeta$, and whose image is $\Cal O_{\bold P^1}(-b)$. In
  this situation, $\tau'$ corresponds to an extension of $i$ to
  $\widetilde Y$, which factors through an embedding of a ribbon
  $\widehat Y$ with conormal bundle $\Cal O_{\bold P^1}(-b)$. In
  particular, $\widetilde Y$ would admit a projection onto $\widehat
  Y$. Now, to prove 1, assume that $\widetilde Y$ does not admit a
  projection onto such $\widehat Y$. Then any preimage $\tau$ of $\zeta$
  by $\delta$ must correspond to a matrix whose first row is
  nonzero. Adding as before a suitable element from the kernel of
  $\delta$ we can actually assume that $\tau$ is a surjective homomorphism of $\text{Hom}(\Cal I/\Cal I^2, \Cal E)$, so $i$ can be lifted to an embedding of $\widetilde Y$ into $\bold P^N$.

\smallskip

Finally, to prove 2, assume that $\widetilde Y$ admits a projection $p$ onto a ribbon $\widehat Y$ with conormal bundle $\Cal O_{\bold P^1}(-b)$.
As before $\widetilde Y$ corresponds to an element $\zeta \in \text{Ext}^1(\Omega_{\bold P^1}, \SE)$. Since Ext$^1(\Omega_{\bold P^N}|_{\bold P^1}, \Cal O_{\bold P^1}(-b))=0$, $i$ can be extended to a morphism $\hat i$ from $\widehat Y$.
Composing  $\hat i$ with $p$, we get also an extension of $i$ to $\widetilde Y$. This extension correspond to an element $\tau \in  \text{Hom}(\Cal I/\Cal I^2, \Cal E)$ whose corresponding matrix has zero entries in its first row. Now any other $\tau'$ mapping to $\zeta$ will differ from $\tau$ by an element in the kernel of $\delta$. Arguing as in the previous paragraph we see that the matrix of such $\tau'$ will have also zeros in its first row. This means that all extensions of $i$ will factor through ribbons whose conormal is a subsheaf of $\Cal O_{\bold P^1}(-b)$. This completes the proof of 2.
\end{proof}

\bigskip

As noted in Remark~\ref{linearly.normal}, the embeddings in Theorem~\ref{embedding} are almost never linearly normal.
However, as corollary of Theorem~\ref{embedding}, we obtain an effective bound for linearly normal embeddings of ropes in projective space:

\begin{cor}\label{higher.embedding}
Let $\WY$ be a rope on $Y=\bold P^1$ with arithmetic genus $g$, Maroni invariant $n$ and conormal bundle $\SE=\Cal O_{\bold P^1}(-a) \oplus \Cal O_{\bold P^1}(-b)$, $a \geq b$.
Let $N_0=\max\{4,a-1\}$ (i.e., $N_0=\max\{4, \frac{g+n}{2}\}$). Let $Y  \overset{i} \hookrightarrow \bold P^N \subset \bold P^M $ be an embedding of $Y$ as a rational normal curve in $\bold P^N$, where $M=N+h^0(\SE\otimes \SO_{\PP^1}(N))$.
If $N \geq N_0$, then there are embeddings
$\WY  \overset{\tilde j} \hookrightarrow  \PP^M$ extending $i$ such
that $\tilde j$ embeds $\WY$ as a linearly normal subscheme of $\bold
P^M$.
\end{cor}

\begin{remark}\label{nonspecial}
{\rm Under the hypotheses of Corollary~\ref{higher.embedding} we have $H^1(\SE \otimes \SO_{\PP^1}(N))=0$.}
\end{remark}

\begin{Proof} (of Corollary~\ref{higher.embedding})
Consider one of the embeddings $\WY \overset{\tilde i} \hookrightarrow \bold P^N$ obtained in Theorem~\ref{embedding}.
Consider also the sequence
\begin{equation*}
0 \longrightarrow H^0(\SE \otimes \SO_{\PP^1}(N)) \longrightarrow H^0(\tilde i^*\SO_{\bold P^N}(1)) \longrightarrow  H^0( \SO_{\PP^1}(N)) \longrightarrow 0,
\end{equation*}
which is right exact by Remark~\ref{nonspecial}. Then the complete linear series of $\tilde i^* \SO_{\bold P^N}(1)$ induces a morphism $\WY \overset{\tilde j} \longrightarrow \bold P^M$, which is an embedding because $\tilde i$ also is.
\end{Proof}\\

\section{Smoothing of multiplicity $3$ ropes on $\bold P^1$}
\label{triplesmooth}
In this section we apply the results of Sections ~\ref{generalsmooth} and ~\ref{embeddings} to prove, in a very precise way,
that ropes on $\PP^1$, of multiplicity $3$ and nonnegative arithmetic genus are smoothable.
Let $\widetilde Y$  be a rope of multiplicity $3$ on $\bold P^1$ and let $\Cal E=\SO_{\PP^1}(-a)\oplus \SO_{\PP^1}(-b)$ be its conormal bundle.
If we want to apply Theorem
~\ref{embsmoothing} to $\widetilde Y$, we need a smooth, irreducible triple cover $X \overset{\pi} \longrightarrow \bold P^1$, with trace zero module $\Cal E$. This imposes conditions on $a$ and $b$, that we summarize in the following well-known result (see ~\cite{Miranda}):

\begin{theorem}\label{triplecovers}
Let $\SE=\SO_{\PP^1}(-a)\oplus \SO_{\PP^1}(-b)$ be a vector bundle of rank two on $\bold P^1$. There exists a smooth, irreducible cover $X \overset{\pi} \longrightarrow \bold P^1$ with $\pi_*\Cal O_X=\Cal O_{\bold P^1} \oplus \Cal E$ if and only if $a>0$, $b>0$, $a \leq 2b$ and $b \leq 2a$. In this case, the general triple cover with trace zero module $\SE$ is smooth and irreducible.
\end{theorem}

Theorems~\ref{embsmoothing} and ~\ref{triplecovers} allow us to say at
once that ropes of genus $g \geq 0$, embedded in projective space,
with a conormal bundle $\SE$ satisfying the above conditions, can be
smoothed.
These are ropes whose Maroni invariant $n$ is small compared to the
genus (precisely, $n \leq \frac {g+2}{3}$). However, this restriction
can be circumvented by degenerating more balanced ropes (to which a
smooth triple cover can be associated) to less balanced ropes (to
which a smooth triple cover cannot be associated). This is done in the
following definition and propositions:

\begin{definition}\label{relropedef}
Let  $\SY$ and $U$ be reduced and connected $\bold k$--schemes,
let $\SY \overset{p} \longrightarrow U$ be a smooth proper
morphism with connected fibers and let $\mathbb{E}$ be a locally free sheaf
on $\SY$. A rope on $\SY$ relative over $U$, with conormal bundle
$\mathbb{E}$,  is a scheme over $U$,  $\widetilde \SY
\overset{\tilde p} \longrightarrow U$, where $\widetilde{\SY}$ is a
rope on $\SY$ with conormal bundle
$\mathbb{E}$ and $p$ is the composition of the closed immersion $\SY
\hookrightarrow \widetilde \SY$ followed by $\tilde p$.
\end{definition}

\begin{proposition}\label{relrope}
Let $\SY \overset{p} \longrightarrow U$,  $\widetilde \SY
\overset{\tilde p} \longrightarrow U$ and $\mathbb{E}$ be as in Definition
~\ref{relropedef}.

\begin{enumerate}
\item The ropes on $\SY$ relative over $U$ with conormal bundle
  $\mathbb{E}$, up to equivalence of extensions, are classified by the
  space $\mathrm{Ext}^1_{\SY}(\Omega_{\SY/U}, \mathbb{E})$.

\item If $U$ is affine, then there is a natural isomorphism between $\mathrm{Ext}^1_{\SY}(\Omega_{\SY/U}, \mathbb{E})$ and $H^0(R^1p_*(\Omega^*_{\SY/U}
  \otimes \mathbb{E}))$.

\item The morphism $\widetilde\SY \overset{\tilde p} \longrightarrow U$ is
  a proper and flat family of ropes $\widetilde\SY_u$ on $\SY_u$ with
  conormal bundle $\mathbb{E}_u= \mathbb{E} \otimes \Cal O_{\SY_u}$,
  $u \in U$.

\medskip

\noindent Let $\Cal Z$ be a smooth, proper scheme over $U$, let $\SY \overset{i} \hookrightarrow \Cal Z$ be a $U$--morphism which
is a closed immersion of smooth varieties and let $\Cal I$ be the
ideal sheaf of $i(\SY)$ in $\Cal Z$.

\medskip

\item The $U$--morphisms $\tilde i$ extending $i$
to $\widetilde \SY$ are classified by $\mathrm{Hom}_{\SY}(\Cal I/\Cal I^2,
\mathbb{E})$.

\item If $(\widetilde \SY, \tilde i)$ corresponds to  $\tau \in
\mathrm{Hom}_{\SY}(\Cal I/\Cal I^2,  \mathbb{E})$, then
$\tilde i$ is a
closed immersion if and only if $\tau$ is
a surjective homomorphism. In this case, $\tilde i$ is a family of
closed immersions of ropes.

\item If $(\widetilde \SY, \tilde i)$ corresponds to  $\tau \in
\mathrm{Hom}_{\SY}(\Cal I/\Cal I^2, \mathbb{E})$ and
$\widetilde \SY$ corresponds to $\varepsilon \in  \mathrm{Ext}^1_{\SY}(\Omega_{\SY/U}, \mathbb{E})$, then $\varepsilon$ is the image of $\tau$
  by the connecting homomorphism.
\end{enumerate}
\end{proposition}

\begin{proof}
1) is a relative version of Theorem ~\ref{paropes}, (1) and the proof
given in ~\cite[Theorem 1.2]{BE}, can be
adapted to the relative setting without problem.

Now we prove 2). We consider the Leray Spectral Sequence
$E_2^{ij}=H^i(R^jp_*\Cal M)$, where $\Cal M=\Omega_{\SY/U}^* \otimes
\mathbb{E}$.  The sequence abuts
to $H^l(\SY, \Cal M)$ and gives the five term exact sequence

\begin{equation}\label{Leray}
0 \longrightarrow H^1(p_*\Cal M) \longrightarrow H^1(\Cal M)
\overset{\rho} \longrightarrow  H^0(R^1p_*\Cal M) \longrightarrow
H^2(p_*\Cal M) \longrightarrow H^2(\Cal M).
\end{equation}

Since $U$ is affine and $p$ is smooth, $\rho$ becomes an isomorphism

$$\mathrm{Ext}^1_{\SY}(\Omega_{\SY/U},\mathbb{E})
\overset{\sim}\longrightarrow H^0(R^1p_*(\Omega_{\SY/U}^* \otimes
\mathbb{E})).$$

Now we prove 3). The morphism $\tilde p$ is proper because $p$ is and
it is also flat because $p$ is and $\mathbb{E}$ is locally free. On the
other hand, by ~(\ref{Leray}), the restriction of an element $\varepsilon \in \mathrm{Ext}^1_{\SY}(\Omega_{\SY/U},\mathbb{E})$ to
$\Cal Y_u$ yields an element $\varepsilon_u \in H^1(\Omega^*_{\SY_u} \otimes \mathbb{E}_u)=
\mathrm{Ext}^1(\Omega_{\SY_u}, \mathbb{E}_u)$. This implies that the
fiber of $\widetilde \SY$ over $u$ is the rope on $\SY_u$, with
conormal bundle $\mathbb{E}_u$ corresponding to $\varepsilon_u$.

Finally (4), (5) and (6) are relative versions of Theorem
~\ref{paropes}, (2), (3) and (4). The proof given in \cite
[Proposition 2.1]{Gonzalez03} translates word by word to the relative
setting.

\end{proof}

\begin{proposition}\label{ropedegen}
Let $a, b$ be integers, $a \geq b+2$. Let $\Cal E'=\Cal O_{\bold
  P^1}(-a+1) \oplus \Cal O_{\bold P^1}(-b-1)$ and let $\Cal E=\Cal
O_{\bold P^1}(-a) \oplus \Cal O_{\bold P^1}(-b)$. Let $\WY$ be a rope
on $Y=\bold P^1$ with conormal bundle $\Cal E$.
\begin{enumerate}
\item
 There exists a flat family $\Cal {\widetilde Y}$ over a neighborhood  $T$ of $0 \in \bold A^1$
such that $\widetilde{\Cal Y_0}=\WY$ and $\widetilde{\Cal Y_t}$ is a rope of conormal bundle $\Cal E'$ if $t \neq 0$.

\item
Let
$N \geq \max\{4,a-1\}$.
Let $Y \overset{i} \hookrightarrow \bold P^N$ be an embedding of $Y$
as a rational normal curve of degree $N$ in $\bold P^N$. Then there
are embeddings $\WY
\overset{\tilde i} \hookrightarrow \bold P^N$ of $\WY$
extending $i$. For such an embedding $\tilde i$,
$T$ and the family $\widetilde{\Cal Y}$ over $T$
described in (1) can be chosen in such a way that $\tilde i$ extends to a family of embeddings $\widetilde{\Cal Y}  \hookrightarrow \bold P_T^N$.
\end{enumerate}
\end{proposition}

\begin{proof}
It is well known (see for instance ~\cite{Harris}) that $\Cal E'$ can
be degenerated to $\Cal E$. More precisely, consider the family $\Cal
Y= \bold P^1 \times \bold A^1 \overset{p_1}\longrightarrow \bold
A^1$. There is a vector bundle $\mathbb{E}$ of rank $2$ on $\Cal Y$ such that the
restriction of $\mathbb{E}$ to $\bold P^1 \times \{0\}$ is $\Cal E$ and the restriction
$\mathbb{E}_t$ of $\mathbb{E}$ to $\bold P^1 \times \{t\}$ is $\Cal E'$ if $t \neq 0$. The
rope $\WY$ corresponds to an element $\zeta \in
\text{Ext}^1(\Omega_{\bold P^1}, \Cal E)=\text{H}^1(\Omega_{\bold
  P^1}^* \otimes
\Cal E)$.
A rope with conormal bundle $\Cal E'$ corresponds to an element of
 $\text{Ext}^1(\Omega_{\bold P^1}, \Cal E')=
 \text{H}^1(\Omega_{\bold P^1}^* \otimes \Cal E')$. Let $p_2$ be the
projection of $\bold P^1 \times \bold A^1$ onto $\bold P^1$. By
Proposition ~\ref{relrope}, a
section of $R^1{p_1}_*( p_2^*(\Omega_{\bold P^1}^*) \otimes \mathbb{E})$
yields a family $\Cal {\widetilde Y}$, flat over $\bold A^1$, such
that $\Cal {\widetilde Y}_0$ is a rope with conormal bundle $\mathbb{E}_0=\SE$
and $\Cal {\widetilde Y}_t$ is a rope with conormal bundle $\mathbb{E}_t=\SE'$ if $t
\neq 0$.  Since the fibers of $\Cal Y \overset{p_1}\longrightarrow
\bold A^1$ have dimension $1$, $R^1{p_1}_*(p_2^*(\Omega_{\bold P^1}^*)
\otimes \mathbb{E})(t)=H^1(\Omega_{\bold
  P^1}^* \otimes \mathbb{E}_t)$ for all $t \in \bold A^1$ (see
~\cite[II.5, Corollary
3]{Mumford}).  Then we can lift $\zeta \in
\text{H}^1(\Omega_{\bold P^1}^* \otimes \Cal E)$ to a section of the
restriction of $R^1{p_1}_*( p_2^*(\Omega_{\bold
  P^1}^*) \otimes \mathbb{E} )$ to a suitable neighborhood $T$ of $0$.  This proves (1).

Now we prove (2). The existence of embeddings $\tilde i$
follows by Theorem~\ref{embedding}. Let $\Cal I$ be the ideal sheaf of
$i(Y)$ in $\bold P^N$. The embedding $\tilde i$ corresponds to a surjective
homomorphism $\tau_0 \in \text{Hom}(\Cal I/\Cal I^2, \mathbb{E}_0)$.
On the other hand, the conditions satisfied by $N$
imply that $h^0({\Cal I/\Cal I^2}^* \otimes \mathbb{E}_t)$ is constant for all $t \in \bold A^1$.
Then $\tau_0$ can be lifted to a local section $\tau$ of ${p_1}_*(p_2^*({\Cal I/\Cal I^2}^*) \otimes \mathbb{E})$.
Since $p_1$ is proper, there is an open neighborhood $T$
of $0 \in \bold A^1$ such that $\tau$ is surjective on $p_1^{-1} \, T$.
Then Proposition~\ref{relrope} implies (2).
\end{proof}

Now we are ready to prove the main theorem of this section:

\begin{theorem}\label{main.rope.theorem}
Let $\WY$ be a  rope on $Y=\bold P^1$ with arithmetic genus $g \geq 0$, Maroni invariant $n$ and conormal bundle $\SE=\Cal O_{\bold P^1}(-a) \oplus \Cal O_{\bold P^1}(-b)$, $a \geq b$.
Let $N_0=\max\{4,a-1\}$ (i.e., $N_0=\max\{4, \frac {g+n}{2}\}$). Let $Y  \overset{i} \hookrightarrow \bold P^N \subset \bold P^M $ be an embedding of $Y$ as a rational normal curve in $\bold P^N$, where $M=N+h^0(\SE\otimes \SO_{\PP^1}(N))$.
If $N \geq N_0$, then

\begin{enumerate}
\item there exist embeddings $\WY \overset{\tilde i} \hookrightarrow
  \bold P^N$ extending $i$; for each such $\tilde i$, there exist an
  affine curve $T$ and a flat family $\widetilde{\Cal Y}$ of subschemes of $\bold
  P^N$ such that $\widetilde{\Cal Y_0}=\tilde i(\WY)$ and
  $\widetilde{\Cal Y_t}$ is a smooth,
  irreducible curve if $t \neq 0$;

\item there exist embeddings
$\WY  \overset{\tilde j} \hookrightarrow  \PP^M$ extending $i$ such
that $\tilde j$ embeds $\WY$ as a linearly normal subscheme of $\bold
P^M$; for each such $\tilde j$, there exist an affine curve $T$ and a
flat family $\widetilde{\Cal Y}$ of subschemes of $\bold P^M$ such that
$\widetilde{\Cal
Y_0}=\tilde j(\WY)$ and $\widetilde{\Cal Y_t}$ is a smooth, irreducible curve if $t
\neq 0$.
\end{enumerate}
\medskip
In particular, all  ropes on $\bold P^1$ with nonnegative arithmetic genus are smoothable.
\end{theorem}

\begin{proof} {\it Case 1.} We assume first that $n \leq
  \frac{g+2}{3}$.
This implies, together with Theorem~\ref{triplecovers}, the existence of a smooth irreducible cover $X \overset{\pi} \longrightarrow Y$ with $\pi_*\Cal O_X=\Cal O_{\bold P^1} \oplus \Cal E$.  On the other hand, the bound on $N$ in terms of $n$ and $g$ implies $\Cal E \otimes \Cal O_Y(1)$ is non--special.
Then the result follows from Theorem~\ref{embsmoothing}, Theorem~\ref{embedding} and Corollary~\ref{higher.embedding}.

\smallskip

{\it Case 2.} We deal now with the remaining ropes. Proposition~\ref{ropedegen}, (1) tells that $\WY$ sits into a flat family $\widetilde{\Cal Y}$ over some affine curve $T$, whose other members are more balanced ropes $\WY'$. More precisely, the ropes $\WY'$ have all the same Maroni invariant $n-2$.
Now, since $N \geq \max\{4,a-1\}$, by Proposition~\ref{ropedegen}, (2), there exists a family of embeddings $\widetilde{\Cal Y} \hookrightarrow \bold P^N_T$.
Arguing as in Corollary ~\ref{higher.embedding}, we have also a family
$\widetilde{\Cal Y} \hookrightarrow \bold P^M_T$ of morphisms
induced by complete
linear series.
Let now $\Cal H$ be the Hilbert scheme of subschemes of $\bold P^N$ of
dimension $1$, degree $3N$ and arithmetic genus $g$. The Hilbert point
of $\WY$ lies in the closure of the locus parametrizing ropes with
Maroni invariant $n-2$, and eventually, in the closure of ropes with
Maroni invariant satisfying the conditions of Case 1. Since the latter
lie in the closure parametrizing smooth curves, so does $\WY$.  Hence
$\tilde i(\WY) \subset \bold P^N$ is smoothable. The same argument
shows that $\tilde j(\WY) \subset \bold P^M$ is smoothable.
\end{proof}

We end by remarking explicitly, as we did in Section~\ref{generalsmooth}, that there are many different smoothings that we can exhibit for a given rope $\WY$:

\begin{remark}\label{manysmooth.triple}
{\rm Given a rope $\WY$ of multiplicity $3$ on $\bold P^1$,  there are many different smoothings of $\WY$
obtained according to the results of this article.
\begin{enumerate}
\item First, as seen in Section~\ref{embeddings}, there are many different embeddings of $\WY$ in projective space.
\item Second, in the case the conormal of $\WY$ satisfies the conditions of Theorem~\ref{triplecovers}, then there are many possible smooth covers $\pi$ to choose (see ~\cite{Miranda}).
\item Third, once $\pi$ is chosen and $\varphi$ is therefore fixed,
  the element $\tau \in \mathrm{Hom}(\SN_{Y,Z}^*,\SE)$ corresponding
  $(\tilde Y, \tilde i)$ may have, in most cases, many different
  liftings to $H^0(\Cal N_\varphi)$. Liftings $\mu$ of $\tau$ to
  $H^0(\Cal N_\varphi)$ correspond to first--order infinitesimal
  deformations $\widetilde \varphi$ of $\varphi$, such that
  $(\mathrm{im}\,\widetilde \varphi)_0=\WY$. For details of this see
  ~\cite[Section 3]{Gonzalez03}.
\item Fourth, there are many different ways of extending $\widetilde \varphi$ to a family of morphisms $\Phi$, as the proof of Theorem~\ref{coversmoothing} shows.
\end{enumerate}}
\end{remark}

\providecommand{\bysame}{\leavevmode\hbox to3em{\hrulefill}\thinspace}
\providecommand{\MR}{\relax\ifhmode\unskip\space\fi MR }
\providecommand{\MRhref}[2]{%
  \href{http://www.ams.org/mathscinet-getitem?mr=#1}{#2}
}
\providecommand{\href}[2]{#2}

\end{document}